\documentclass{gtmon_a}  
\pdfoutput=1
\usepackage{amscd}
\numberwithin{equation}{section}

\proceedingstitle{Groups, homotopy and configuration spaces (Tokyo
  2005)}
\conferencestart{5 July 2005}
\conferenceend{11 July 2005}
\conferencename{Groups, homotopy and configuration spaces,
                in honour of Fred Cohen's 60th birthday}
\conferencelocation{University of Tokyo, Japan}

\editor{Norio Iwase}
\givenname{Norio}
\surname{Iwase}

\editor{Toshitake Kohno}
\givenname{Toshitake}
\surname{Kohno}

\editor{Ran Levi}
\givenname{Ran}
\surname{Levi}

\editor{Dai Tamaki}
\givenname{Dai}
\surname{Tamaki}

\editor{Jie Wu}
\givenname{Jie}
\surname{Wu}

\title[The homology of spaces of polynomials]{The homology of spaces of polynomials\\
with roots of bounded multiplicity} 

\author{Yasuhiko Kamiyama}
\givenname{Yasuhiko}
\surname{Kamiyama}
\address{Department of Mathematics, University of the Ryukyus\\ 
Nishihara-Cho, Okinawa 903-0213, Japan}
\email{kamiyama@sci.u-ryukyu.ac.jp}
\urladdr{}

\keyword{polynomial}
\keyword{loop space}
\keyword{integral homology}
\subject{primary}{msc2000}{55P35}
\subject{secondary}{msc2000}{55R20}
\subject{secondary}{msc2000}{58D15}

\arxivreference{}

\volumenumber{13}
\issuenumber{}
\publicationyear{2008}
\papernumber{12}
\startpage{281}
\endpage{291}
\doi{}
\MR{}
\Zbl{}

\received{5 June 2006}
\revised{20 February 2007}
\accepted{27 February 2007}
\published{25 February 2008}
\publishedonline{25 February 2008}
\proposed{}
\seconded{}
\corresponding{}
\version{}

\makeatletter
\def\cnewtheorem#1[#2]#3{\newtheorem{#1}{#3}[section]
\expandafter\let\csname c@#1\endcsname\c@thm}

\theoremstyle{plain}
  \newtheorem{thm}{Theorem}[section]
  \cnewtheorem{lem}[thm]{Lemma}
  \cnewtheorem{cor}[thm]{Corollary}
  \cnewtheorem{prop}[thm]{Proposition}
  
\theoremstyle{definition}
  \cnewtheorem{defn}[thm]{Definition}
  \cnewtheorem{exmp}[thm]{Example}
  \cnewtheorem{rem}[thm]{Remark}
  \cnewtheorem{rems}[thm]{Remarks}

\makeatother  
\makeautorefname{defn}{Definition}

\begin{document}

\begin{htmlabstract}
Let P<sup>l</sup><sub>k,n</sub> be the space consisting of monic complex
polynomials f(z) of degree k and such that the number of n&ndash;fold roots
of f(z) is at most l. In this paper, we determine the integral homology groups
of P<sup>l</sup><sub>k,n</sub>.
\end{htmlabstract}

\begin{abstract}  
Let $P_{k, n}^l$ be the space consisting of monic complex
polynomials $f(z)$ of degree $k$ and such that the number of $n$--fold roots 
of $f(z)$ is at most $l$. In this paper, we determine the integral homology groups 
of $P_{k, n}^l$. 
\end{abstract}

\begin{asciiabstract}
Let P_{k, n}^l be the space consisting of monic complex polynomials
f(z) of degree k and such that the number of n-fold roots of f(z) is
at most l. In this paper, we determine the integral homology groups of
P_{k, n}^l.
\end{asciiabstract}

\maketitle

\section{Introduction} \label{sec1}
In \cite{A}, Arnol'd studied a space $P_{k, n}^l$ 
consisting of monic complex
polynomials $f(z)$ of degree $k$ and such that the number of $n$--fold roots of $f(z)$ 
is at most $l$. In particular, he calculated the first five nontrivial 
integral homology groups of $P_{k, n}^l$. The purpose of this paper is, 
using another approach, to 
determine $H_\ast (P_{k, n}^l; \mathbb{Z})$ completely.
\par
Let $C_k (\mathbb{C})$ denote the configuration space of unordered $k$--tuples 
of distinct points in $\mathbb{C}$. The study of the topology of $C_k (\mathbb{C})$ 
originated in \cite{A}. For that purpose, Arnol'd performed an induction for 
$P_{k, n}^l$ with making $k$ larger and $l$ smaller while $n$ being fixed. 
Then one obtains information on $P_{k, n}^l$ for all $k$, $n$ and $l$. 
In particular, setting $n=2$ and $l=0$, 
we obtain information on $C_k (\mathbb{C})$. 
(Strictly speaking, Arnol'd considered the complement $S^{2k}- P_{k, n}^l$ 
instead of $P_{k, n}^l$.) 
\par
Using this induction, Arnol'd calculated the first five nontrivial 
integral homology groups of $P_{k, n}^l$. (See \fullref{3.1} for $n=2$.) But because of
problems involved in the induction, it seems difficult to calculate
further homology groups. Then we naturally encounter the following problem: 
how to determine $H_\ast (P_{k, n}^l; \mathbb{Z})$. 
\par
The purpose of this paper is to give an answer to the problem. Our main results
will be stated in \fullref{sec3}. (See Theorems \ref{3.3} and \ref{3.7}.) Here we summarize how the groups
$H_\ast (P_{k, n}^l; \mathbb{Z})$ are determined. 
\begin{thm}\label{1.1}
Let $J^l (S^{2n-2})$ be the $l$-th stage of the James construction 
which builds $\Omega S^{2n-1}$, and let $W^l (S^{2n-2})$ be the homotopy theoretic fiber 
of the inclusion $J^l (S^{2n-2}) \hookrightarrow \Omega S^{2n-1}$. Then:
\begin{enumerate}
\renewcommand{\labelenumi}{\rm(\roman{enumi})}
\item 
\begin{enumerate}
\renewcommand{\labelenumi}{\rm(\alph{enumi})}
\item
The homomorphism
\begin{equation*}
H_\ast (P_{k,n}^l; \mathbb{Z}) \to H_\ast (P_{k+1,n}^l; \mathbb{Z})
\end{equation*}
which is induced from the natural inclusion 
$P_{k,n}^l \hookrightarrow P_{k+1,n}^l$ is a monomorphism onto a direct summand.
\item
There is a stable homotopy equivalence
\begin{equation*}
P_{\infty, n}^l \underset{s}{\simeq} W^l (S^{2n-2}).
\end{equation*}
\end{enumerate}
\item The homology groups $H_\ast (W^l (S^{2n-2}); \mathbb{Z})$ are determined. 
In particular, all higher $p$--torsions are determined for all primes $p$. 
\item
For each $x \in H_\ast (W^l (S^{2n-2}); \mathbb{Z})$, the least $k$ 
such that $x$ is contained in $H_\ast (P_{k,n}^l; \mathbb{Z})$ is determined.
\end{enumerate}
\end{thm}

\begin{rem}\label{1.2}
For $l=0$, \fullref{1.1} is already well-known.  First, about
\fullref{1.1} (i) (a), the inclusion $P_{k,n}^0 \hookrightarrow P_{k+1,n}^0$,
which is called a stabilization map, was constructed by Guest,
Kozlowski and Yamaguchi in \cite{GKY1,GKY2}. Moreover, the induced
homomorphism $H_\ast (P_{k,n}^0; \mathbb{Z}) \to H_\ast (P_{k+1,n}^0;
\mathbb{Z})$ was studied in \cite{GKY2}.   Second, about \fullref{1.1} (i)
(b) and (iii) for $l=0$, Guest, Kozlowski and Yamaguchi \cite{GKY1}
and independently Kallel \cite{Kal1} established a more precise
result. (See \fullref{2.2}.)
\end{rem}

Finally, we note that the homology groups $H_\ast (C_k (\mathbb{C});
\mathbb{Z}/p)$ were determined later, using other approaches, by Fuks
for $p=2$ \cite{F} and by F~Cohen for odd primes $p$ \cite{C}. F~Cohen
also determined Steenrod operations.  \par This paper is organized as
follows. In \fullref{sec2} we summarize previous results on $P_{k,
n}^l$ which imply \fullref{1.1} (i). In \fullref{sec3} we first recall
Arnol'd's results in \fullref{3.1}. Our main result for $n=2$ is
\fullref{3.3}, which generalizes \fullref{3.1}.  \fullref{3.7} is a
generalization of \fullref{3.3} for general $n$.

\par
The author would like to thank the referee for invaluable suggestions.

\section{Previous results}\label{sec2}
As in \fullref{sec1}, we set
\begin{align*}
P_{k, n}^l = \{ f(z) &:\text{$f(z)$ is a monic complex polynomial 
of degree $k$}\\
&\text{and such that the number of $n$--fold roots of $f(z)$ is at most $l$}\}. 
\end{align*}
Since $P_{k,n}^l = \mathbb{C}^k$ for $k< n (l+1)$, 
we can assume that $k \geq n(l+1)$. 
\par
On the other hand, let $\text{\rm Rat}_k (\mathbb{C} P^{n-1})$ denote the 
space of based holomorphic maps of degree $k$ from the Riemannian sphere $S^2$ 
to the complex projective space $\mathbb{C} P^{n-1}$. The basepoint  condition 
we assume is that $f (\infty)=[1, \dots, 1]$. Such holomorphic maps are given 
by rational functions:
\begin{align*}
\text{\rm Rat}_k (\mathbb{C} P^{n-1}) = \{ (p_1 (z),& \dots, p_n (z)): 
\text{each $p_i(z)$ is a monic degree--$k$ polynomial}\\
&\text{and such that there are no roots common to all $p_i (z)$}\}.
\end{align*}
The study of the topology of $\text{\rm Rat}_k (\mathbb{C} P^{n-1})$ originated 
in Segal's paper \cite{S}, where it is proved that the natural inclusion 
$\text{\rm Rat}_k (\mathbb{C} P^{n-1}) \hookrightarrow \Omega_k^2 \mathbb{C} P^{n-1}
\simeq \Omega^2 S^{2n-1}$ is a homotopy equivalence up to dimension $k (2n-3)$. 
\par
Later, F~Cohen et al determined the stable homotopy type of 
$\text{\rm Rat}_k (\mathbb{C} P^{n-1})$ as follows: 
\begin{thm}{\rm\cite{CCMM1,CCMM2}}\label{2.1}\qua
Let
\begin{equation*}
\Omega^2 S^{2n-1} \underset{s}{\simeq} \underset{1 \leq j}{\bigvee} D_j (S^{2n-3})
\end{equation*}
be Snaith's stable splitting. Then there is a stable homotopy equivalence
\begin{equation*}
\text{\rm Rat}_k (\mathbb{C} P^{n-1}) \underset{s}{\simeq}
\underset{j=1}{\overset{k}{\bigvee}}
D_j (S^{2n-3}).
\end{equation*}
\end{thm}
In particular, combining \fullref{2.1} for $n=2$ with the stable splitting of 
$C_k (\mathbb{C})$ (Brown and Peterson \cite{BP}), we have 
\begin{equation}\label{2-1}
C_k (\mathbb{C)} \underset{s}{\simeq}
\text{\rm Rat}_{\left[ \frac{k}{2} \right]} (\mathbb{C} P^1).
\end{equation}
Guest, Kozlowski and Yamaguchi and independently Kallel generalized \eqref{2-1} as follows: 
\begin{thm}{\rm\cite{GKY1,Kal1}}\label{2.2}\qua
For $n \geq 3$, there is a homotopy equivalence
\begin{equation*}
P_{k, n}^0 \simeq \text{\rm Rat}_{\left[ \frac{k}{n} \right]} (\mathbb{C} P^{n-1}).
\end{equation*}
\end{thm}
\begin{rems}\label{2.3}
(i)\qua It is proved by Guest, Kozlowski and Yamaguchi in \cite{GKY2}
that the (modified) jet map $P_{k,n}^0 \to \text{\rm Rat}_k
(\mathbb{C} P^{n-1})$ defined by
\begin{equation*}
f(z) \mapsto (f(z), f(z)+f'(z), \dots, f(z)+f^{(n-1)}(z))
\end{equation*}
is a homotopy equivalence up to dimension $(2n-3) \left[\frac{k}{n} \right]$ if
$n \geq 3$, and a homology equivalence up to dimension 
$(2n-3) \left[ \frac{k}{n} \right]$ if $n=2$.
\par
(ii)\qua Kallel \cite{Kal2} generalized $P_{k,n}^0$ as follows: 
let $F^d (\mathbb{R}^m, k)$ be the space of 
ordered $k$--tuples of vectors in $\mathbb{R}^m$ 
so that no vector occurs more than $d$ times in the $k$--tuple. 
We set $C^d (\mathbb{R}^m, k) = 
F^d (\mathbb{R}^m, k)/\Sigma_k$. 
Then $C^1 (\mathbb{R}^m, k)$ is the usual configuration space and 
$C^{n-1}(\mathbb{R}^2,k) \cong P_{k,n}^0$. 
Recall that using 
$F^1 (\mathbb{R}^m, k)$, May, Milgram and Segal constructed 
a combinatorial model for $\Omega^m \Sigma^m X$, 
where $X$ is a connected CW--complex. 
Using $F^d (\mathbb{R}^m,k)$, 
Kallel \cite{Kal2} generalized the model for general $d$. He also
considered the case when $X$ is disconnected. In particular, 
setting $m=2$, $d=n-1$ and $X= S^0$ in his result, he recovered
the homotopy and homology equivalences
$P_{\infty, n}^0 \simeq \Omega^2 S^{2n-1}$ for $n \geq 3$ and $n=2$, 
respectively. (See \fullref{2.2} and \eqref{2-1} for these equivalences.) 
\par
(iii)\qua For $n \geq 2$, a stable homotopy equivalence
\begin{equation}\label{2-2}
P_{k, n}^0 \underset{s}{\simeq} \text{\rm Rat}_{\left[ \frac{k}{n} \right]} 
(\mathbb{C} P^{n-1})
\end{equation}
was proved by Vassiliev in \cite{V}. \fullref{2.2} is a stronger version of \eqref{2-2} for $n \geq 3$.
\end{rems}

We consider generalizations of Theorems \ref{2.1} and \ref{2.2}. We set 
\begin{align*}
X_{k, n}^l =  \{ (p_1 (z),& \dots, p_n (z)): 
\text{each $p_i(z)$ is a monic degree--$k$ polynomial}\\
&\text{and such that there are at most $l$ roots common to all $p_i (z)$} \}.
\end{align*}
\begin{thm}{\rm(Kamiyama \cite{Kam1})}\label{2.4}\qua
Let $J^l (S^{2n-2})$ denote the $l$-th stage of the James 
construction which builds $\Omega S^{2n-1}$, and let $W^l (S^{2n-2})$ be the
homotopy theoretic fiber of the inclusion $J^l (S^{2n-2}) \hookrightarrow 
\Omega S^{2n-1}$. Let
\begin{equation*}
W^l (S^{2n-2}) \underset{s}{\simeq} \underset{1 \leq j}{\bigvee} D_j \xi^l (S^{2n-2})
\end{equation*}
be a generalization of Snaith's stable splitting. (See Wong \cite{W}
and Kamiyama \cite{Kam1}.)
Then, there is a stable homotopy equivalence
\begin{equation*}
X_{k, n}^l \underset{s}{\simeq}
\underset{j=1}{\overset{k}{\bigvee}}
D_j \xi^l (S^{2n-2}).
\end{equation*}
\end{thm}
\begin{thm}[Kamiyama \cite{Kam2}]\label{2.5}
For $l \geq 1$ and $n \geq 2$, there is a homotopy equivalence
\begin{equation*}
P_{k, n}^l \simeq X_{\left[ \frac{k}{n} \right], n}^l.
\end{equation*}
\end{thm}
Note that \fullref{1.1} (i) are consequences of Theorems 2.4 and 2.5. 

\section{The main results}\label{sec3}
In order to simplify notation, we first consider the case $n=2$, 
which is of particular interest to us. 
Since $P_{k, 2}^l = \mathbb{C}^k$ for $k < 2l+2$, we assume that
$k \geq 2l+2$. 
\par
Arnol'd proved the following:
\begin{thm}{\rm\cite{A}}\label{3.1}\ 

\begin{enumerate}
\renewcommand{\labelenumi}{\rm(\roman{enumi})}
\item 
For $1 \leq j \leq 2l$, we have $H_j (P_{k, 2}^l; \mathbb{Z}) =0$.
\item
For $2l+1 \leq	j \leq 2l+5$, the groups $H_j (P_{k, 2}^l; \mathbb{Z})$ 
are cyclic and the orders are given by the following table.
\end{enumerate}

\begin{table}[ht!] \small
\caption{The orders of the groups $H_j (P_{k, 2}^l; \mathbb{Z}) \;
(2l+1 \leq j \leq 2l+5)$}
\begin{center}
\begin{tabular}{cccccc}
\noalign{\hrule height0.8pt} 
\hfil $k$ $\backslash$ $j$ & $2l+1$ & $2l+2$ & $2l+3$ & $2l+4$ & $2l+5$\\
\hline
$2l+2, 2l+3$ & $\infty$ & $0$ & $0$ & $0$ & $0$\\
$2l+4, 2l+5$ & $\infty$ & $l+2$ & $0$ & $0$ & $0$\\
$2l+6, 2l+7$ & $\infty$ & $l+2$ & $2/(l+1)$ & 
$(l+3)/2$ & $0$\\
$2l+8, 2l+9$ & $\infty$ & $l+2$ & $2/(l+1)$ 
& $((l+3)/2)(2/(l+1))$ & $3/(l+1)$\\
$2l+10, 2l+11$ & $\infty$ & $l+2$ & $2/(l+1)$ 
& $((l+3)/2)(2/(l+1))$ & $6/(l+1)$\\
$\vdots$ & $\vdots$ & $\vdots$ & $\vdots$ & $\vdots$ & $\vdots$\\
$\infty$ & $\infty$ & $l+2$ & $2/(l+1)$ 
& $((l+3)/2)(2/(l+1))$ & $6/(l+1)$\\
\noalign{\hrule height0.8pt} 
\end{tabular}
\end{center}
\end{table}
Here we introduce the notation 
\begin{equation*}
a/b = \frac{a}{\gcd (a, b)},
\end{equation*}
where $\gcd (a, b)$ is the greatest common divisor 
of the integers $a$ and $b$.
\end{thm}
\par
In order to state our main results, we prepare some notation.
\begin{defn}\label{3.2}
Let $p$ be a prime. 
\begin{enumerate}
\renewcommand{\labelenumi}{\rm(\roman{enumi})}
\item 
We write $l$ as $l= p^m q$ such that 
\begin{equation*}
q= \sum_{\nu=0}^N a_\nu p^\nu,
\end{equation*}
where $0 \leq a_\nu \leq p-1$ and $a_N \not= 0$, $a_0 \not=0$.
\item
For $q$ in (i), we consider terms of the form
\begin{equation*}
(p-1) \sum_{\nu=j}^i p^\nu.
\end{equation*}
We take such terms as large as possible, whence we have 
$a_\nu = p-1 \; ( j \leq \nu \leq i)$ and $a_{i+1} \not= p-1$, $a_{j-1} \not= p-1$.
Assume that all possible pairs $(i, j)$ for $q$ are given by 
\begin{equation*}
(i_\alpha, j_\alpha), \quad 1 \leq \alpha \leq r,
\end{equation*}
where we arrange them as $j_\alpha \geq i_{\alpha+1}+2$.
\item For $1 \leq \alpha \leq r$, we set
\begin{equation*}
u_\alpha = \sum_{\nu=i_\alpha +1}^N a_\nu p^\nu.
\end{equation*}
\item
We set 
\begin{equation*}
d_\alpha = 2 \left( p^m u_\alpha+ p^{m+i_\alpha+1}-1 \right).
\end{equation*}
\item
We set
\begin{equation*}
\mu_\alpha = i_\alpha- j_\alpha+2.
\end{equation*}
\end{enumerate}
\end{defn}
Our main result for $n=2$ is then:
\begin{thm}\label{3.3} Let $p$ be a prime. Then all higher $p$--torsions in 
$H_\ast (W^l (S^2); \mathbb{Z})$ are given as follows.
\begin{enumerate}
\renewcommand{\labelenumi}{\rm(\roman{enumi})}
\item 
If $m \geq 1$, then
\begin{enumerate}
\renewcommand{\labelenumi}{\rm(\alph{enumi})}
\item
For $1 \leq \alpha \leq r$, 
$H_{d_\alpha} (W^l (S^2); \mathbb{Z})$
contains  
$\mathbb{Z}/p^{\mu_\alpha}$ as a direct summand. 
\item
For each $\alpha$, the least $k$ such that the higher $p$--torsion in (a) 
appears as a direct summand in $H_{d_\alpha} (P_{k, 2}^l; \mathbb{Z})$ is
\begin{equation*}
k= d_\alpha +2.
\end{equation*}
\end{enumerate}
\item
If $m=0$, then we omit the case $\alpha=r$ from (i). 
\end{enumerate}
\end{thm}
\begin{rem}\label{3.4}
We can determine all $p$--torsions of order exactly $p$ in 
$H_\ast(P_{k, 2}^l; \mathbb{Z})$ from the following facts:
all $p$--torsions in $H_\ast (W^l (S^2); \mathbb{Z})$ of order exactly $p$ 
are determined from the Bockstein operation on $H_\ast (W^l (S^2); \mathbb{Z}/p)$, 
and $H_\ast (P_{k, 2}^l; \mathbb{Z}/p)$ is a subspace of 
$H_\ast (W^l (S^2); \mathbb{Z}/p)$ (see \fullref{3.6}).
Hence using \fullref{3.3}, we know the groups $H_\ast (P_{k, 2}^l; \mathbb{Z})$ completely.
\end{rem}
\begin{exmp}\label{3.5}
We consider the case 
\begin{equation*}
l= p^m (p-1) \left( \sum_{\nu=j_1}^{i_1} p^\nu+ \sum_{\nu=j_2}^{i_2} p^\nu \right).
\end{equation*}
\begin{enumerate}
\renewcommand{\labelenumi}{\rm(\roman{enumi})}
\item 
If $m \geq 1$, then there are $2$ higher $p$--torsions:
\begin{enumerate}
\item For $k \geq 2p^{m+i_1+1}$, 
\begin{equation*}
H_{2(p^{m+i_1+1}-1)} (P_{k,2}^l; \mathbb{Z})
\end{equation*}
contains $\mathbb{Z}/p^{i_1-j_1+2}$ as a direct summand. 
\item 
For $k \geq 2p^m (p^{i_1+1}-p^{j_1}+p^{i_2+1})$, 
\begin{equation*}
H_{2p^m (p^{i_1+1}-p^{j_1}+p^{i_2+1})-2} (P_{k,2}^l; \mathbb{Z})
\end{equation*}
contains $\mathbb{Z}/p^{i_2-j_2+2}$ as a direct summand. 
\end{enumerate}
\item If $m=0$, then we omit the case (b) from (i).
\end{enumerate}
\end{exmp}
\begin{proof}[Proof of \fullref{3.3}]
(i)\qua In order to prove (a), we determine $H_\ast (W^l (S^2); \mathbb{Z})$ 
by the following $2$ steps. 
\begin{enumerate}
\item
Using the structure of $H_\ast (W^l (S^2); \mathbb{Z}/p)$, we determine the
homological dimensions which have higher $p$--torsions.
\item
Using the cohomology Serre spectral sequence for a fibration 
with coefficients in $\mathbb{Z}_{(p)}$, 
we determine the higher $p$--torsions. 
\end{enumerate}
\par
(1)\qua The structure of $H_\ast (W^l (S^2); \mathbb{Z}/p)$ was determined in \cite{Kam1}
from the mod $p$ Serre spectral sequence for the fibration
\begin{equation*}
\Omega^2 S^3 \to W^l (S^2) \to J^l (S^2).
\end{equation*}
Let $x \in H_2 (J^l (S^2); \mathbb{Z}/p)$ and 
$\iota \in H_1 (\Omega^2 S^3; \mathbb{Z}/p)$ be the generators
and we write $Q_1^t = Q_1 \cdots Q_1$ ($= t$-times $Q_1$). 
In $H_\ast (W^l (S^2); \mathbb{Z}/p)$, the cases that the Bockstein operation is not
clear are given as follows:
\begin{equation}\label{3-1}
x^{p^m u_\alpha} \otimes Q_1^{m+i_\alpha+1} (\iota) \to 
x^{p^m v_\alpha} \otimes \beta Q_1^{m+j_\alpha} (\iota), 
\quad 1 \leq \alpha \leq r,
\end{equation}
where we set 
\begin{equation*}
v_\alpha = \sum_{\nu= j_\alpha}^N a_\nu p^\nu.
\end{equation*}
(Note that by \fullref{3.2}, we have $v_\alpha = u_\alpha+ 
\sum_{\nu=j_\alpha}^{i_\alpha} (p-1) p^\nu$. Note also that $v_\alpha= u_{\alpha+1}$
for $p=2$.) Since
\begin{equation*}
\deg (x^{p^m v_\alpha} \otimes \beta Q_1^{m+j_\alpha} (\iota))= d_\alpha,
\end{equation*}
there is a higher $p$--torsion in $H_{d_\alpha} (W^l (S^2);\mathbb{Z})$. 
This completes (1). 
\par
(2)\qua Consider the following homotopy commutative diagram:
\begin{equation*}
\begin{CD}
W^l (S^2) @>>> \tilde{J}^l (S^2) @>>>\Omega S^3 \langle 3 \rangle\\
\Vert @. @VVV @VVV\\
W^l (S^2) @>>> J^l (S^2) @>>> \Omega S^3\\
@VVV @VVV @VVV\\
\ast @>>> K(\mathbb{Z}, 2) @>>> K(\mathbb{Z}, 2)
\end{CD}
\end{equation*}
where $\tilde{J}^l (S^2)$ and $\Omega S^3 \langle 3 \rangle$ are the homotopy 
theoretic fibers of the second and third columns respectively. 
Then the first row is a fibration and we consider 
the cohomology Serre spectral sequence for the fibration 
with coefficients in $\mathbb{Z}_{(p)}$. 
Note that $H^{d_\alpha+1} (W^l (S^2); \mathbb{Z}_{(p)})$ is determined 
if we calculate the cokernels of the differentials 
\begin{equation}\label{3-2}
d: E^{2ps, d_\alpha-2ps+1} \to E^{d_\alpha+2, 0}
\end{equation}
for all possible $s \geq 1$. Since $H_q (W^l (S^2); \mathbb{Z}_{(p)}) =0$ for
$q \leq 2l$, we have the following restriction on $s$:
$d_\alpha -2ps+1 \geq 2l+1$, that is, 
\begin{equation}\label{3-3}
p^{m+j_\alpha}-1- \sum_{\nu=0}^{m+j_\alpha-2} b_\nu p^\nu \geq ps,
\end{equation}
where $0 \leq b_\nu \leq p-1$. 
\par
Let $y_{2ps} \in H^{2ps} (\Omega S^3 \langle 3 \rangle; \mathbb{Z}_{(p)})$ be
a generator. Then a generator of $E^{2ps, d_\alpha-2ps+1}$ is mapped by $d$ 
in \eqref{3-2} to $y_{2ps}y_{d_\alpha-2ps+2}$. It is easy to see that 
\begin{equation}\label{3-4}
y_{2ps}y_{d_\alpha-2ps+2} = {p^m u_\alpha+p^{m+i_\alpha+1} \choose ps}y_{d_\alpha+2}.
\end{equation}
Consider the $p$--power component of the prime decomposition of the binomial 
coefficient in \eqref{3-4}. Using \eqref{3-3}, we see that the component is smallest when 
$ps = t p^{m+j_\alpha-1} \, (1 \leq t \leq p-1)$ such that the $p$--power is 
$p^{i_\alpha-j_\alpha+2}$. Hence $$H^{d_\alpha+1} (W^l (S^2); \mathbb{Z}_{(p)}) =
 \mathbb{Z}/p^{\mu_\alpha}$$ and \fullref{3.3} (i) (a) follows. 
\par
For \fullref{3.3} (i) (b), we have the following:
\begin{prop}\label{3.6}
In $H_\ast (P_{k,2}^l; \mathbb{Z}/p)$, we define the weights of the homology classes 
$x$ and $\iota$ (see \eqref{3-1}) to be $2$. Then $H_\ast (P_{k,2}^l; \mathbb{Z}/p)$ is 
isomorphic to the subspace of $H^\ast (W^l (S^2); \mathbb{Z}/p)$ spanned by 
monomials of weight $\leq k$. 
\end{prop}
\begin{proof} 
The proposition is an easy consequence of Theorems 2.4 and 2.5. 
Note that it is reasonable to define the weights of $x$
and $\iota$ to be $2$ by the following reason: we have $H_\ast (W^l (S^2); \mathbb{Q}) = 
\bigwedge (x^l \otimes \iota)$. Since $P_{2l+2,2}^l \simeq S^{2l+1}$, the weight of
$x^l \otimes \iota$ must be $2l+2$. 
\end{proof}
Since the weight of 
$x^{p^m v_\alpha} \otimes \beta Q_1^{m+j_\alpha} (\iota)$ in \eqref{3-1} is $d_\alpha+2$, 
\fullref{3.3} (i) (b) follows. 
\par
(ii)\qua For $m=0$ and $\alpha=r$, the left-hand side of \eqref{3-1} is the mod $p$ reduction
of the generator of $H_{2l+1} (W^l (S^2); \mathbb{Z}) = \mathbb{Z}$ and the right-hand
side is $0$. Hence we must omit this case from (i). 
This completes the proof of \fullref{3.3}.
\end{proof}

Finally we generalize \fullref{3.3} for general $n$. 
\begin{thm}\label{3.7} We keep the notation of \fullref{3.2} except that we 
generalize $d_\alpha$ in (iv) as
\begin{equation*}
d_{n, \alpha} = 2 (n-1)p^m (u_\alpha+ p^{i_\alpha+1})-2.
\end{equation*}
Then: 
\begin{enumerate}
\item
\fullref{3.3} (i) (a) is generalized to the assertion that 
$H_{d_{n, \alpha}} (W^l (S^{2n-2}); \mathbb{Z})$ contains 
$\mathbb{Z}/p^{\mu_\alpha}$ as a direct summand. 
\item
About \fullref{3.3} (i) (b), the least $k$ 
such that the higher $p$--torsion in the above (1) 
appears as a direct summand in $H_{d_{n, \alpha}} (P_{k, n}^l; \mathbb{Z})$ is
\begin{equation*}
k= \frac{n (d_{n,\alpha} +2)}{2 (n-1)}.
\end{equation*}
\item
\fullref{3.3} (ii) holds under these modifications. 
\end{enumerate}
\end{thm}
\begin{proof} 
About $x$ and $\iota$ in \eqref{3-1}, we generalize that $x \in 
H_{2n-2} (J^l (S^{2n-2}); \mathbb{Z}/p)$ and $\iota
\in H_{2n-3} (\Omega^2 S^{2n-1}; \mathbb{Z}/p)$ such that 
the weights of these elements are $n$. 
\fullref{3.7} is clear from this.
\end{proof}

\bibliographystyle{gtart}
\bibliography{link}

\end{document}